\RequirePackage{fix-cm}
\documentclass[smallextended]{svjour3} 
\smartqed 
\usepackage{graphicx}

\usepackage{color}

\usepackage{amsmath}
\usepackage{amssymb}
\usepackage[pagebackref]{hyperref} 
\usepackage{framed} 
\usepackage{empheq} 

\newcommand{\pathprefix}[1]{#1}

\newcommand{\wh}{\widehat}
\newcommand{\ol}{\overline}

\newcommand{\R}{\mathbb{R}}
\newcommand{\eps}{\varepsilon}

\newtheorem{prop}{Proposition}

\begin{document}

\title{Sparse solutions of optimal control via Newton method for under-determined systems
\thanks{This work was supported by Russian Science Foundation (project 16-11-10015).}
}

\titlerunning{Sparse control} 

\author{Boris Polyak \and
 Andrey Tremba 
}


\institute{B. Polyak \at Institute for Control Sciences, Profsoyuznaya 65, 117997 Moscow, Russia,
\email{boris@ipu.ru} 
\and
A. Tremba \at Institute for Control Sciences, Profsoyuznaya 65, 117997 Moscow, Russia,
\email{atremba@ipu.ru} 
}

\date{Received: date / Accepted: date}

\maketitle

\begin{abstract}
We focus on finding sparse and least-$\ell_1$-norm solutions for unconstrained nonlinear optimal control problems. Such optimization problems are non-convex and non-smooth, nevertheless recent versions of Newton method for under-determined equations can be applied successively for such problems.
\end{abstract}

\vskip 32pt

\section{Introduction}
The simplest optimal control problem has the form
\begin{equation}\label{opt}
\begin{array}{c}
\displaystyle
\min \int_{0}^{T} F(x,u,t)dt\\[3mm]
\displaystyle
\dot{x}=\phi (x,u,t), \; x(0)=a, \;x(T)=b\in \R^m, \; u(t) \in \R^q
\end{array}
\end{equation}
where there are no constraints on state $x(t)$ and control $u(t)$. There are numerous classical methods for solving the problem such as gradient method, dynamic programming, maximum principle etc. see \cite{BN,Bellman,Evt,LM,Leitman,Leondes,Polak}. However these methods are oriented on smooth functionals. Our problems of interest are sparse  controls, i.e. solutions with small number of non-zero components. Such solutions naturally arise in $L_1$ control, where the performance criteria are $L_1$ norms of control (or where $L_1$ regularizer is exploited). The $L_1$ norm is $\int_{t=0}^T \|u(t)\|_1 dt$, it involves vector $\ell_1$-norm (Manhattan, $1$-norm) $\|u(t)\|_1 = \sum_{i=1}^q |u_i(t)|$. Examples are minimal fuel control for satellite trajectory optimization, where $L_1$ term corresponds to the fuel consumption \cite{Athans,TK}, {Section 5.5} or so called ``maximal hands-off'' control \cite{NQN}. If we take finite-dimensional approximation of such problems we arrive to mathematical programming problems with nonlinear (and hence non-convex) equality-type constraints and non-smooth objective function. The $L_1$ norm transforms to $\ell_1$-norm of the extended vector of (discretized) control. There are no standard algorithms for such hard problems. Of course situation is much simpler if the state equations are linear, then linear programming technique can be applied, see \cite{Rao} and references therein. 

The contributions of the paper are twofold. First, we exhibit that recently developed by the authors Newton-like method for under-determined equations \cite{Polyak-Tremba,Polyak-Tremba-OPTIMA}  can be highly effective for solving nonlinear optimal control problems. That is, we consider right-hand side condition $x(T)=b$ as $m$ nonlinear equations with variables $u$ (the number of variables $n$ depends on discretization and in general is much larger $m$). Above mentioned method iteratively solves linearized equation, finding solution with minimal norm. Second, we apply the method to find sparse solutions. This is due the flexibility of the proposed method to various norms. By choosing $L_1$ norm we can find sparse controls. Recent paper \cite{PMY} also exploits methods from \cite{Polyak-Tremba,Polyak-Tremba-OPTIMA} for optimal control, but with $L_2$ norm of approximations, thus the solutions are not sparse.

Theoretically, the solution of the problem \eqref{opt} with $L_1$ terms may contain impulses, e.g. a combination of delta-functions. However, differential equation is discretized in practice, and the following problem with embedded difference equation is considered
\begin{gather}\label{eq:discrete-objective}
\displaystyle
\min F(\{x\}, \{u\})\\[2mm]
\label{eq:discrete-equation}
\begin{array}{c}
\displaystyle
x[j+1] = f_j(x[j], u[j]), \; j = 0, ..., N-1, \\[2mm]
\displaystyle
x[0]=a, \;x[N] = b\in \R^m,
\end{array}
\end{gather}
with functions $f_j : \R^{m + q} \rightarrow \R^m$. The functions $f_j$ are assumed to be differentiable on both arguments, while $F$ can be non-differentiable.
 Notice that the dimension of the variables in the problem becomes very large: $q N$, with $N \approx T/h \gg 1$, $h$ being (small) discretization step-size.

If the objective equals $L_1$ norm of control in continuous-time setup, it becomes simple $\ell_1$ norm in the discretized problem statement:
\begin{equation} \label{eq:obj-l1}
F(\{x\}, \{u\}) = \sum_{j=0}^{N-1} \|u[j]\|_1.
\end{equation}
The function summarizes absolute values of all components for all time instances.
In algorithmic approach, we deal with the discretized equations only.

Another objective function 
 is the \emph{total number} of non-zero components, counted through all time instants:
\begin{equation} \label{eq:obj-total-sparsity}
F(\{x\}, \{u\}) = \text{total \# of non-zero components in } \{u[j]\}_{j=0,...,N-1}.
\end{equation}

We aim to solve problem \eqref{eq:discrete-objective}-\eqref{eq:discrete-equation} with one of the above objective functions. Running ahead, the optimization cannot be done exactly, but an approximate algorithm is proposed, based on the specially fitted version of Newton method. 
We discuss its applicability and properties. 
It appears that this approach can manage sparsity constraint \eqref{eq:obj-total-sparsity} directly.

The structure of the paper is as follows. In Section 2 we describe classical Newton method and its versions for under-determined equations. Section 3 contains main results of the paper and addresses applications of the methods to sparse control. There are several versions of the algorithms to get sparse solutions. Results of numerical simulation are provided in Section 4.

\section{Constraint Equation and Newton method}
In this section equation \eqref{eq:discrete-equation} and its solution are considered.
It is convenient to vectorize the whole control sequence $\{u[j]\}, j = 0, ...,N-1$, stacking individual vectors $u[j] \in \R^q$ into a big one.
Then the control sequence $\{u[k]\}$ is represented by the single vector variable $u = (u[0]^T, u[1]^T, \ldots, u[N-1]^T)^T \in \R^{q N}$. Let's denote its dimension as $n = qN$.

Equations \eqref{eq:discrete-equation} can be presented in the form of the single equation
\begin{equation} \label{eq:P_u}
P(u) \doteq x[N](u) - b = 0 \in \R^m, \; u \in \R^n.
\end{equation}
The function $P(u)$ is easily calculated by applying the recursive relation of \eqref{eq:discrete-equation}:
$$
x[j+1] = f_j(x[j], u[j]), \; j = 0, ..., N-1, \; x[0] = a.
$$
Equation \eqref{eq:P_u} is strongly under-determined, as $m \ll n = q N$.
If the functions $f_j$ are differentiable, so is the function $P$. Its derivative is easily obtained by the chain rule applied to $P(u) = f_{N-1}(f_{N-2}(..., u[N-2]), u[N-1])$, see an example below in Section~\ref{sec:opt-control}.

\subsection{Regular Newton method and its versions}

One of the most popular generic algorithms for solving smooth nonlinear equations is Newton method, which uses the idea of the linearization of $P(u)$ at each iteration \cite{DS,Kelley,OR}:
$$
u^{k+1} = u^k - \gamma_k (P'(u^k))^{-1} P(u^k), \; k = 0, 1, ...
$$
The formulae includes both pure ($\gamma_k \equiv 1$) and damped ($\gamma_k \leq 1$) variants of Newton method. 

It is known that if Newton method converges, then its convergence rate is quadratic $\|P(u^{k+1})\| \leq c_1 \|P(u^{k})\|^2$.
Famous Newton-Kantorovich theorem impose semi-local conditions, ensuring that pure Newton method converges \cite{Kant82}. 

The theorem operates with four entities, matching four assumptions:
\begin{enumerate}
\item[$\mathbf{A}$.]
Let $\mu_0$ be a constant, describing non-degeneracy $P'$ at a \emph{single point} $u^0$. 
\item[$\mathbf{B}$.]
Let $L$ be a constant, describing variability of $P'$ \emph{around} $u^0$.
\item[$\mathbf{C}$.] 
Let $\rho$ be the radius of the abovementioned ``around $u^0$'' area, where the constant $L$ is valid.
\item[$\mathbf{D}$.]
And let $s$ be the size of initial residual, e.g. $\|P(u^0)\|$.
\end{enumerate}
We intentionally avoid unneccessary formalization of the constants here, to clarify them with respect to under-determined case later. In a particular case, $\mu_0$ is the least singular value of $P'(u^0)$ and $L$ is the Lipschitz constant of the derivative.
It is sufficient for $P$ to be differentiable around $u^0$, but for simplicity of statements we assume that $P$ is differentiable everywhere.

We also omit the exact formulation of the Newton-Kantorovich theorem, but its virtue is the following: 

If $\mathbf{A}, \mathbf{B}, \mathbf{C}, \mathbf{D}$ hold, and the following conditions 
\begin{equation} \label{eq:newton-kant-cond}
h \doteq \frac{L}{\mu_0^2}s < \frac{1}{2}, \; \text{ and }
\frac{1 - \sqrt{1-2h}}{h} \frac{s}{\mu_0} \leq \rho,
\end{equation}
hold true, then pure Newton method converges with quadratic rate \cite{Kant82}. There are multiple results on rigorous convergence conditions, cf. surveys \cite{Polyak2006,Yam}. 

Another theorem, by Mysovskikh \cite{Kant82}, 
exploits the following assumption $\mathbf{A'}$ instead of $\mathbf{A}$:
\begin{enumerate}
\item[$\mathbf{A'}$.]
Let $\mu$ be a constant, uniformly describing non-degeneracy of $P'$ \emph{around} point $u^0$.
\end{enumerate}
The typical result \cite[Corollary 1]{Polyak-1964} reads as: 
If $\mathbf{A'}, \mathbf{B}, \mathbf{C}, \mathbf{D}$ hold, and inequalities
\begin{equation} \label{eq:newton-mys-cond}
h \doteq \frac{L}{\mu^2} s < 2, \; \text{ and } 
\frac{2\mu}{L} H_0\left( \frac{h}{2} \right) < \rho,
\end{equation}
are satisfied, then pure Newton method 
converges to the solution $u^*$ with quadratic convergence rate. 
The monotonically increasing function $H_0(\delta) = \sum_{\ell=0}^{\infty} \delta^{(2^\ell)}$ is an infinite sum of double exponents.

In both cases, from conditions \eqref{eq:newton-kant-cond},\eqref{eq:newton-mys-cond} it follows, that convergence radius $s = \|P(u^0)\|$ is limited for fixed constants $L$ and $\mu$ (or $\mu_0$). It means that the convergence of pure Newton method is always local.

Damped Newton method is a very natural choice for situations when pure Newton method fails. The damping parameter $\gamma_k \leq 1$ allows to make non-unit step-size in Newton direction $-w^k = -(P'(u^k))^{-1} P(u^k)$, and due to the step-size tuning damped Newton method may be converging globally. Algorithms for adjusting $\gamma_k$ can be found in \cite{Burd,OR}.

\subsection{Newton method for under-determined equations}
Regular Newton method assumes the existence of the inverse operator $(P'(u))^{-1}$
near the solution and initial point, so direct application of the method is restricted to the case $n = m$.

The first version of the Newton method for arbitrary $m \neq n$ case has been proposed in \cite{BI}, see also \cite{L_BI}:
\begin{equation} \label{eq:pseudo}
u^{k+1} = u^k - (P'(u^k))^\dagger P(u^k),
\end{equation}
here $A^\dagger$ is the Moore-Penrose pseudo-inverse of a matrix $A$. This includes under-determined case $m < n$. Other versions of the method for under-determined equations can be found in \cite{Hager,Nesterov,Walker} and references therein. We rely on the approach, proposed in our previous works \cite{Polyak-Tremba,Polyak-Tremba-OPTIMA}.
Its general form looks as follows
\begin{equation} \label{eq:newton-w}
\begin{array}{l}
\displaystyle
w^k = \arg \min_{P'(u^k) w = P(u^k)} \|w\|, \\
\displaystyle
u^{k+1} = u^k - \gamma_k w^k.
\end{array}
\end{equation}
This description unifies under-determined and regular case. 
Two key properties of \eqref{eq:newton-w} are:
\begin{enumerate}
\item 
freedom of choosing a norm of $w$ to be minimized, and
\item
freedom of step-size $\gamma_k$ choice.
\end{enumerate}
The first property is completely absent in regular, well-defined systems of equations with $n = m$. For under-determined case, if norm is Euclidean, $w^k = (P'(u^k))^\dagger P(u^k)$, and we arrive to the method with pseudo-inverse operator \eqref{eq:pseudo}.
For our purposes (finding sparse solutions) we use method in form \eqref{eq:newton-w} with $\ell_1$-norm, with zero initial condition $u^0 = 0$. Then the auxiliary problem of finding $w^k$ direction is equivalent to linear programming one.

In \cite{Polyak-Tremba,Polyak-Tremba-OPTIMA} the authors proposed two novel, alternative strategies for choosing damping parameter $\gamma_k$:
\begin{equation} \label{eq:optimal-gamma-muL}
\gamma_k = \min\Big\{1, \; \frac{\mu^2}{L \|P(x^k)\|}\Big\},
\end{equation}
and
\begin{equation} \label{eq:optimal-gamma-L}
\gamma_k = \min\Big\{1, \; \frac{\|P(u^k)\|}{L \|w^k\|^2}\Big\}.
\end{equation}
Whenever needed, these strategies perform damped Newton steps at the beginning, switching to pure Newton method later. 
There are adaptive variants of the strategies, requiring no knowledge of the parameters $\mu, L$ (of $\mathbf{A'}, \mathbf{B}$), see Algorithm in Section~\ref{sec:opt-control}.

\subsection{Convergence Results}
For an under-determined equation $P(u) = 0, \; P : \R^n \rightarrow \R^m$ with $m \leq n$, we are to establish few results for ``sparse-Newton'' method \eqref{eq:newton-w} in the form
\begin{equation} \label{eq:newton-w-l1}
\begin{array}{l}
\displaystyle
w^k = \arg \min_{P'(u^k) w = P(u^k)} \|w\|_1, \\
\displaystyle
u^{k+1} = u^k - \gamma_k w^k.
\end{array}
\end{equation}
with the initial condition $u^0$.
Main topics of interest is the behavior of the objective functions \eqref{eq:obj-l1} and \eqref{eq:obj-total-sparsity}.

Assumptions/constants for the specific case are the following:
\begin{enumerate}
\item[$\mathbf{C}$.] 
$\rho$ is the radius of the $\ell_1$-ball around the initial point: $B_\rho = \{u : \|u - u^0\|_1 \leq \rho \}$.
\item[$\mathbf{A}$.]
$\mu_0>0$ is a constant, satisfying $\|(P'(u^0))^T h\|_1 \geq \mu_0 \|h\|_\infty $ for all $h \in \R^m$. 
\item[$\mathbf{A'}$.]
$\mu>0$ is a constant, satisfying $\|(P'(u))^T h\|_1 \geq \mu \|h\|_\infty $ for all $h \in \R^m$ and all $u \in B_\rho$. 
\item[$\mathbf{B}$.]
Let $L$ be the Lipschitz constant of $P'(u)$: 
$$\|P'(u^a) - P'(u^b)\|_\infty \leq L \|u^a - u^b\|_1, \;\; \forall u^a, u^b \in B_\rho.$$
\item[$\mathbf{D}$.]
$s$ is the $\ell_\infty$-norm of the initial residual, i.e. $s = \|P(u^0)\|_\infty$.
\end{enumerate}

First we consider pure ``sparse-Newton'' method with $\gamma_k \equiv 1$.

If $\mathbf{A}, \mathbf{B}, \mathbf{C}, \mathbf{D}$ hold alongside with \eqref{eq:newton-kant-cond}, then the ``sparse-Newton'' method \eqref{eq:newton-w-l1} converges to a solution $u^*$, and $\|u^0 - u^*\|_1 \leq \frac{\mu_0}{L} \Big(1 - \sqrt{1 - \frac{2L}{\mu_0^2} s} \Big)$.

If $\mathbf{A'}, \mathbf{B}, \mathbf{C}, \mathbf{D}$ hold alongside with \eqref{eq:newton-mys-cond}, then the ``sparse-Newton'' method \eqref{eq:newton-w-l1} converges to a solution $u^*$, and $\|u^0 - u^*\|_1 \leq \frac{\mu}{L} H_0 \big(\frac{L}{2\mu^2} s\big)$.

The results follow from \eqref{eq:newton-kant-cond}, \eqref{eq:newton-mys-cond} by virtue of analysis \cite{Kant82} for specified norms.
In both cases the initial residual $s$ is limited by $\frac{\mu_0^2}{2L}$ or $\frac{2\mu^2}{L}$.

Newton algorithm \eqref{eq:newton-w} with adaptive step-size \eqref{eq:optimal-gamma-muL} extends the limits for $s$, and being potentially unlimited allows global convergence. Convergence conditions extend \eqref{eq:newton-mys-cond}, cf. \cite{Polyak-Tremba} for details.
Particularly, Newton algorithm with step-size \eqref{eq:optimal-gamma-muL} makes not more than
\begin{equation} \label{eq:k_max-definition} 
k_{\max} = \max\Big\{0, \;\left\lceil \frac{2 L}{\mu^2} s \right\rceil - 2 \Big\} 
\end{equation}
damped Newton iterations, followed by pure Newton iterations.
At $k$-th step the following estimates for the rate of convergence hold:

\begin{equation} \label{eq:dist-x_k-x_star}
\|u^{k} - u^*\|_1 \leq \left\{
\begin{array}{ll}
\displaystyle
\frac{\mu}{L} \big( k_{\max} - k + 2 H_0\Big(\frac{\ol{v}}{2}\Big) \big), & k < k_{\max},\\[2mm]
\displaystyle
\frac{2\mu}{L}H_{k - k_{\max}}\Big(\frac{\ol{v}}{2}\Big), & k \geq k_{\max}.
\end{array}
\right.
\end{equation}
{\color{red}where $\ol{v} = \frac{L}{\mu^2} s - \frac{k_{\max}}{2} < 1$.}

If started from $u^0 = 0$, the upper estimate for \eqref{eq:obj-l1} is 
$$\|u^*\|_1 \leq \frac{\mu}{L} \big( k_{\max} + 2 H_0\Big(\frac{\ol{v}}{2}\Big) \big).$$

Thus if the sparse Newton method converges to the solution, the total number of steps to achieve accuracy $\|P(u^k)\|_\infty \leq \eps$ is estimated as 
$$K_\eps = k_{\max} + \left\lceil \log_2\Big( \log_{\ol{v}/2}\Big(\frac{\eps L}{2\mu^2}\Big)\Big)\right \rceil.$$
As the method achieves the given accuracy, sparsity is also present:
\begin{prop}
For $u^0 = 0$ there exists an approximate solution $\wh{u} \doteq u^{K_\eps}$ with not more than $K_{\eps} m$ non-zero components, such that $P(\wh{u}) = \|x[T](\wh{u}) - b\|_\infty \leq \eps$.
\end{prop}

Indeed each of the $w_i$ in Newton step \eqref{eq:newton-w-l1} is a solution of $\ell_1$-optimization problem with $m$ equality constraints. If the constraints are compatible, then there exists a solution $w_i$ with $m$ non-zero components only \cite{Candes-et-al,Rao}, so up to $k$-th step there are not more than $km$ non-zero components in $u^k$.

This proposition connects two terms: Sparsity and $\ell_1$-objective function. The latter is a common convex substitute to the ``number of non-zero components'' criteria, and is commonly used in machine learning, compressed sensing etc. \cite{Candes-et-al}.
We propose using the same heuristic for objective \eqref{eq:obj-total-sparsity}, but iteration-wise. Coupled with fast convergence of Newton method, it results in the total sparsity of the solution.

\section{Application to sparse control}
\label{sec:opt-control}

Consider a non-linear discrete-time dynamic system on a finite interval: 
\begin{equation} \label{eq:iter-linear-on-u}
x[j+1] = f(x[j]) + B u[j], \;\; j= 0, ..., N-1,
\end{equation}
with $x[j] \in \mathbb{R}^m, u[j] \in \mathbb{R}^q$, differentiable 
function $f: \R^m \rightarrow \R^m$ and matrix $B$ of proper dimension $q\times m$.
Here we use bracketed argument to distinguish time dynamic from algorithmic iterations.
For simplicity we restrict ourselves with systems, nonlinear in $x$ and linear in $u$; the extension to general case \eqref{eq:discrete-equation} is obvious.

Given the initial condition $x[0]$ and the terminal point $x[N] = b$, the goal is to find
a control sequence $u = \{u[j]\}$ which would satisfy condition
$x[N] = b$ via ``the least effort control'', i.e. solving optimization problem
\begin{equation} \label{eq:l1-nonlinear-constr}
\|u\|_1 \rightarrow \min_{P(u) = 0}
\end{equation}
with $P(u)= x[N](u) - b$. 
In the case of a non-linear function $f(x)$ the optimization problem \eqref{eq:l1-nonlinear-constr} is generally very hard to solve. This is due to the equality-type constraint $x[N] = b$ describing essentially non-convex set in the space of control variables.

The problem \eqref{eq:l1-nonlinear-constr} has rather interesting property --- its solution appears to be sparse, i.e. contains few non-zero elements.
However, if a generic optimization algorithm is used for solving \eqref{eq:l1-nonlinear-constr} straightforwardly, e.g. via (sub)gradient method with projections on the manifold $\{u : P(u) = 0\}$, then it meets two obstacles. 
First, finding the projection is hard problem itself, and second --- its solution would be non-sparse.

We propose a heuristic algorithm to find an approximate solution of \eqref{eq:l1-nonlinear-constr}.
The idea is to put sparsity property of solution in front. New problem is informally described as
$$
\text{Find sparse solution } u : P(u) = 0.
$$

That is we are to find a solution with small number of non-zero components.
This can be achieved by multiple ways, for example:
\begin{enumerate}
\item 
Fix all except few coordinates of $u$ to zeros. Then try to solve equation $P(u) = 0$ for the free components only.
If all combinations are checked, the optimal one may be found as well.
\item
Run ``sparse-Newton'' algorithm \eqref{eq:newton-w-l1}, which adds few non-zero components into $u$ at each step.
\item
Combine two previous approaches: Run ``sparse-Newton'' algorithm \eqref{eq:newton-w-l1} for one or few steps, then fix some components of $u$ based on these steps and solve the system with respect to the selected components.
\end{enumerate}

First approach is hard to implement, because of its combinatorial nature.
The second and third approaches rely on special ``sparse-Newton'' algorithm \eqref{eq:newton-w-l1}. The idea is to incorporate our desire of sparsity \emph{within} Newton steps. This can be perfectly done by use of $\ell_1$-norm.
In the third approach number of non-zero components in $u$ is controlled directly, 
and it can be applied for treating objective function \eqref{eq:obj-total-sparsity}.

Important part of the algorithm is the choice of step-size $\gamma_k$. As problem's constants are rarely known, policies \eqref{eq:optimal-gamma-muL}, \eqref{eq:optimal-gamma-L} are not applicable directly. Let's describe an algorithm, requiring no a-priori knowledge of the constants. It is based on \eqref{eq:optimal-gamma-muL}.

The algorithm is initialized with scalar parameters $\beta_0 > 0$, $0 < q < 1$, stopping condition $\eps > 0$, counter $k = 0$, zero initial condition $u^0 = 0 \in \R^n$, and number $p_0 = \|P(u^0)\|_\infty$. 

\begin{framed} \begin{minipage}{0.9\linewidth}
\textbf{Adaptive sparse-Newton algorithm}

\begin{enumerate}

\item \label{it:sparse-opt-step}
Solve LP problem
$\displaystyle w^k = \arg \min_{P'(u^k) w = P(u^k)} \|w\|_1.$

\item \label{it:sparse-step-calc}
Evaluate
$
p_{k+1} = \Big\|P\big(u^k - \min\big\{1, \frac{\beta_k}{p_k}\big\} w^k\big)\Big\|_\infty.
$

\item 
If either a)
$
 \beta_k < p_k \text{ and } p_{k+1} < p_k - \frac{\beta_k}{2}, 
$

or b)
$
 \beta_k \geq p_k \text{ and } p_{k+1} < \frac{1}{2 \beta_k} p_k^2
$ 
holds, then go to Step~\ref{it:sparse-ok-step}.
\item
Update $\beta_{k} \leftarrow q \beta_k$ and return to Step~\ref{it:sparse-step-calc} without increasing counter.

\item \label{it:sparse-ok-step}
Take
$
u^{k+1} = u^k - \min\big\{1, \frac{\beta_k}{p_k}\big\} w^k,
$
set $\beta_{k+1} = \beta_k$, increase counter $k \leftarrow k+1$.
\item
Check stopping condition $p_k \leq \eps$. If it holds, return $u^k$ as the solution. Otherwise return to Step~\ref{it:sparse-opt-step}.
\end{enumerate}
\end{minipage}
\end{framed}
The algorithm has similar to Proposition~1 sparsity property: After $k$ steps there is not more than $k m$ non-zero components in $u^k$. If Newton method converges in few iterations, then resulting $u$ has few non-zero components as well.

The convergence of the adaptive sparse-Newton algorithm is stated in terms of $\beta_0, q, \|P(0)\|_\infty$ and properties $\mathbf{A', B, C, D}$ of $P$, cf. \cite{Polyak-Tremba}. 

\subsection{Calculating derivatives}

Calculating function $P(u)$ and derivative $P'(u)$ is an easy task.
First, run iterations \eqref{eq:iter-linear-on-u}, resulting in intermediate $x[j], j=1,..., N$, which
immediately lead to $P(u) = x[N] - b$.
Then we apply chain rule to
\begin{align*}
P(u) & = x[N](u) - b = f(x[N-1]) + B u[N-1] - b= \\
& = f\big(f(x[N-2]) + B u[N-2]\big) + B u[N-1] - b= \ldots =\\
& = f\Big(f\big(f(...) + B u[N-3]\big) + B u[N-2]\Big) + B u[N-1] - b.
\end{align*}
Blocks of partial derivatives $Q_r = P'_{u[r]}$ contain products of the Jacobian matrices of $f(\cdot)$ as
$$
Q_r=\prod^{N-1}_{j=N-r} f'(x[j]) = f'(x[N-r]) f'(x[N-r+1])\cdot...\cdot f'(x[N-1]) B, r = 1, ..., N-1,
$$
with $Q_0 = B$.
The derivative $P'(u) \in \R^{p \times qN}$ is formed as horizontally concatenated matrix 
\begin{align*}
P'(u) & = [Q_{N-1}, Q_{N-2}, ..., Q_{1}, Q_0] = \\ 
& = \left[\prod^{N-1}_{j=1} f'(x[j]) B, ... , f'(x[N-2]) f'(x[N-1]) B, f'(x[N-1]) B, B\right].
\end{align*}

Generic nonlinear discrete system \eqref{eq:discrete-equation}
can be treated similarly. The only difference is a more complex routine for the evaluation of $P'(u)$.

\section{Example}
Dynamic equation \eqref{eq:iter-linear-on-u} often arises due to discretization of a continuous-time dynamic systems with small step-size. 
As simple example we consider a mathematical pendulum with friction. In normalized variables of angle $\phi$ and angular speed $\psi= \dot{\phi}$ pendulum's dynamics can be described as
$$
\ddot{\phi} + \alpha \dot{\phi} + \beta \sin (\phi) = u,
$$
or
$$
\begin{pmatrix}
\dot{\phi} \\
\dot{\psi}
\end{pmatrix} 
= 
\begin{pmatrix}
\psi \\
- \alpha \psi - \beta \sin(\phi)
\end{pmatrix} 
+ 
\begin{pmatrix}
0 \\ 1
\end{pmatrix} 
u
$$
The latter equation has discrete counterpart with $(x_1[j], x_2[j]) \approx (\phi(j h), \psi(j h) )$. It is obtained by Euler forward method with step $h$. Written in the form \eqref{eq:iter-linear-on-u}, function $f$ and matrix $B$ are
$$
f(x) = 
\begin{pmatrix}
x_{1} + h x_{2} \\ 
x_{2} + h (-\alpha x_{2} - \beta \sin (x_{1}))
\end{pmatrix} 
, \;\; 
B = 
\begin{pmatrix}
0 \\ 1
\end{pmatrix}.
$$
For numerical experiments we set $\alpha = 0.3, \beta = 0.9, h = 0.04, N = 160$,
$x[0] = (1, 0.5)^T, b = (0.4, 0)^T$.

At the first experiment we use sparse-Newton method, following the second approach mentioned above (by running algorithm \eqref{eq:newton-w-l1}). It converges in 3 steps, and after the second step $\|P(x^k)\|_\infty \leq 10^{-5}$. 
At last step machine accuracy is achieved.
Trajectories for all control sequences $u^k, k=0,...,3$ are shown as two upper plots on Figure~1. 
The trajectories of $u^2, u^3$ are visually non-distinguishable.
Components of internal variable $w^k$, used in sparse-Newton method, are seen on the next three plots. There are only two non-zero components at each iteration indeed. Notice that the magnitudes of the components rapidly decrease over iterations.
Resulting control $u^3$ is on the bottom plot of Figure~1. It has only 5 non-zero elements with indices $j \in \{99, 153, 154, 158, 159\}$.
The value of the functional $\|u^3\|_1 = 0.6544$.

\begin{figure*}
\includegraphics[width=1.0\textwidth]{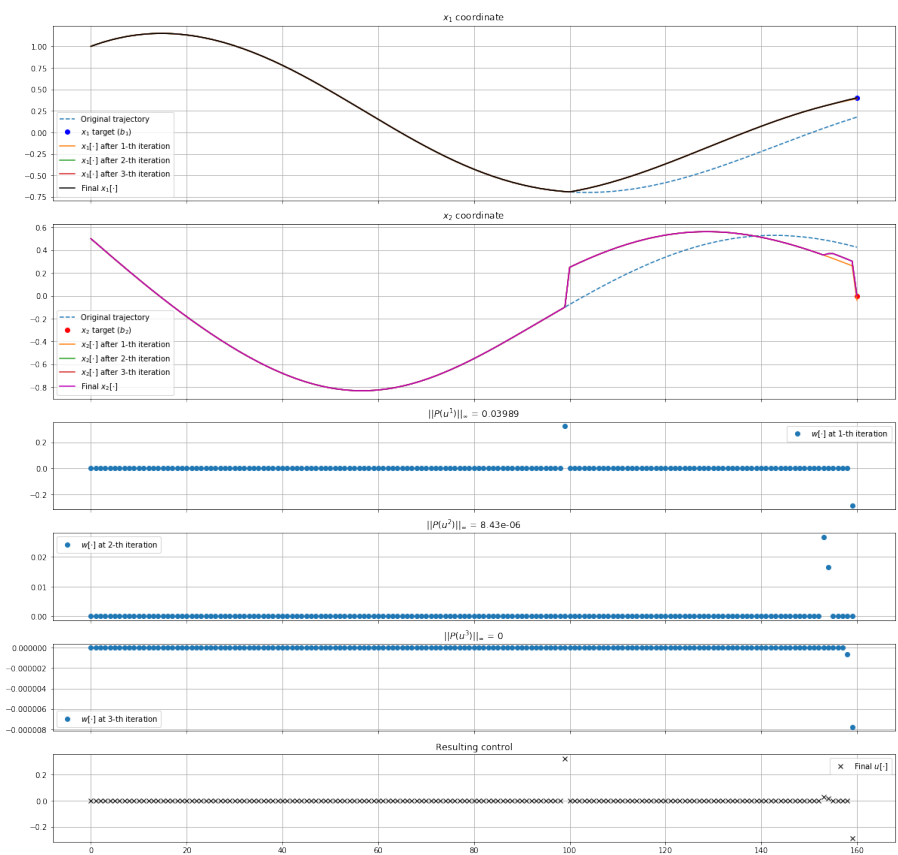}
\caption{Iterations of sparse Newton method \eqref{eq:newton-w-l1}: 5 non-zero components}
\label{fig:2} 
\end{figure*}

For the second experiment we use the third approach with the following tactic: After the first iteration of sparse-Newton method,
the nonzero components $w$ are revealed (indices $j \in \{99, 159\}$). Then we fix all other components of $u$ as zero, and
solve well-defined equation $P(u) = 0$ as $P\big(u=(0, ...,0, u_{99}, 0, ..., 0, u_{159}) \big) = P(u_{99}, u_{159}) = 0$.
The equation is solved by adaptive Newton method with accuracy $10^{-9}$.

The result of the experiment is on Figure~2.
Now the resulting control has functional $\|u\|_1 = 0.6331$, which appears to be less than in the first experiment, with only 2 non-zero components. 

\begin{figure*}
\includegraphics[width=1.0\textwidth]{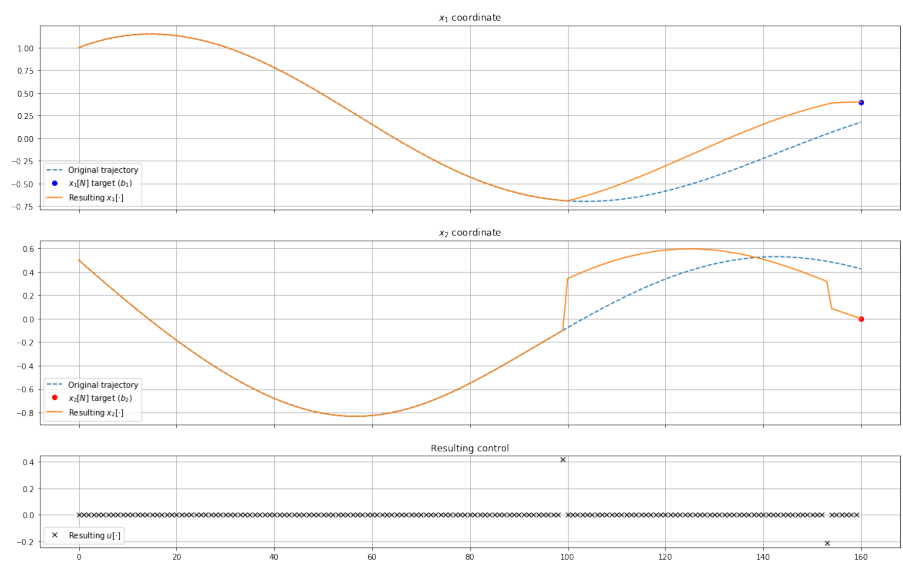}
\caption{Solution with two non-zero coordinates}
\label{fig:2} 
\end{figure*}

\section{Conclusions}
We have proposed the novel approach for finding sparse solutions of boundary-value nonlinear dynamic problems, based on Newton-like method for solving under-determined systems of nonlinear equations. Practical simulation demonstrates high efficiency of the algorithm.

\subsection*{Acknowledgements}

This work was supported by Russian Science Foundation, Project 16-11-10015. The authors thank the anonymous  reviewers for their helpful comments.

\end{document}